\documentclass[reqno,11pt]{amsart}
\usepackage{latexsym}
\usepackage{fullpage}
\usepackage{amssymb}
\usepackage{hyperref}
\def\R{\ensuremath{\mathbb R}}
\def\C{\ensuremath{\mathbb C}}
\def\Z{\ensuremath{\mathbb Z}}
\def\N{\ensuremath{\mathbb N}}
\def\id{\operatorname{id}}
\let\eps=\varepsilon
\def\rk{\operatorname{rk}}
\def\ind{\operatorname{ind}}
\def\tr{\operatorname{tr}}
\def\rk{\operatorname{rk}}
\def\Ad{\operatorname{Ad}}
\def\ch{\operatorname{ch}}
\def\frg{{\mathfrak g}}
\def\frh{{\mathfrak h}}
\def\frp{{\mathfrak p}}
\def\del{\partial}
\def\norm#1{\left\|#1\right\|}
\def\abs#1{\left|#1\right|}
\def\<{\langle}
\def\>{\rangle}
\def\Hom{{\operatorname{Hom}}}
\def\End{{\operatorname{End}}}
\def\SU{\mathrm{SU}}
\def\U{\mathrm U}
\def\SO{\mathrm{SO}}
\def\Spin{\mathrm{Spin}}
\def\so{\mathfrak{so}}
\def\G{\mathrm G}
\def\O{{\mathbb O}}
\def\I{{\mathbb I}}
\def\adsl{\tilde{\mathrm{ad}}}

\numberwithin{equation}{section}
\swapnumbers
\theoremstyle{plain}
\newtheorem{Lemma}[equation]{Lemma}
\newtheorem{Proposition}[equation]{Proposition}

\newtheorem{thm}{Theorem}

\theoremstyle{definition}

\theoremstyle{remark}
\newtheorem{rem}{Remark}

\newtheorem{Example}[equation]{Example}

\begin{document}

\title{Vafa-Witten Estimates for Compact Symmetric Spaces}
\author{S. Goette}
\address{NWF 1 Mathematik, Universit\"at Regensburg,
93040 Regensburg, Germany}
\email{sebastian.goette@mathematik.uni-regensburg.de}
\thanks{Supported in part by DFG special programme
``Global Differential Geometry''}
\begin{abstract}
We give an optimal upper bound for the first eigenvalue
of the untwisted Dirac operator on a compact symmetric space~$G/H$
with~$\rk G-\rk H\le 1$ with respect to arbitrary Riemannian metrics.
We also prove a rigidity statement.
\end{abstract}

\maketitle

Herzlich gave an optimal upper bound for the lowest eigenvalue
of the Dirac operator on spheres with arbitrary Riemannian
metrics in~\cite{herz}
using a method developed by Vafa and Witten in~\cite{VW}.
More precisely, he proved that for every metric~$\bar g$ on~$S^n$
that is pointwise larger than the round metric~$g$,
the first eigenvalue~$\lambda_1(\bar D^2)$ of the Dirac operator
with respect to~$\bar g$ is not larger than the first
Dirac eigenvalue~$\lambda_1(D^2)$ of the round sphere.

Herzlich asked if there are other Riemannian manifolds with optimal
Vafa-Witten bounds,
in particular if the Fubini-Study metric on~$\C P^{2m-1}$ has this property.
In the present note we give positive answers to both questions
by generalising Herzlich's results to symmetric spaces~$G/H$ of compact type,
where~$\rk G-\rk H\le 1$.
In particular,
we  improve a recent estimate by Davaux and Min-Oo
for complex projective spaces in~\cite{DM}, see Example~\ref{CPnEx} below.

\begin{thm}\label{thm1}
  Let~$M=G/H$ be a simply connected symmetric space of compact type
  with~$\rk G-\rk H\le 1$
  and assume that~$M$ is $G$-spin.
  Let~$g$ be a symmetric metric, and let~$D$ denote the corresponding Dirac
  operator on~$M$.
  If~$\bar g$ is another metric with~$\bar g\ge g$ on~$TM$
  and~$\bar D$ is the corresponding Dirac operator,
  then
	$$\lambda_1(\bar D^2)\le\lambda_1(D^2)\;.$$
  In the case of equality,
  we have~$\bar g=g$.
\end{thm}

For an arbitrary Riemannian metric~$\bar g$
such that~$c^2\bar g\ge g$ for some suitable positive constant~$c^2$,
the theorem implies
	$$\lambda_1(\bar D^2)\le c^2\,\lambda_1(D^2)\;.$$

We combine the methods of~\cite{herz} and~\cite{DM}
with related estimates in~\cite{GS}.
In particular,
we compare~$\bar D$ to an operator~$\bar D_1$ with nonvanishing kernel
acting on the same sections as~$\bar D_0=\bar D\otimes\id_{\R^k}$.
We use Parthasarathy's formula to compute~$\lambda_1(D^2)$,
and we exhibit a similar formula to estimate~$\norm{\bar D_1-\bar D}$.
Both formulas give the same value for~$g=\bar g$.

To prove that~$\bar D_1$ has a kernel,
we use the invariance of the Fredholm index if~$\rk H=\rk G$.
If~$\rk H=\rk G-1$ we use the invariance of the mod-2-index as in~\cite{GS}.
Unfortunately,
both approaches fail if~$\rk G-\rk H\ge 2$.
Note that in~\cite{herz}, a spectral flow argument was
used instead in the case~$\rk H=\rk G-1$.

In~\cite{baum}, Baum applied the Vafa-Witten approach to Lipschitz
maps~$f$ of high degree from a closed Riemannian spin manfifold
of dimension~$2n$ to~$S^{2n}$.
We extend her result to Lipschitz maps of high $\hat A$-degree
from higher dimensional closed Riemannian spin manifolds to~$S^{2n}$.
Recall that if~$N^n$ and~$M^m$ are closed oriented manifolds,
$[N]$ is the fundamental class of~$N$
and~$\omega$ is a generator of~$H^m(M;\Z)$,
then the $\hat A$-degree is defined as
	$$\deg_{\hat A}f=\bigl(\hat A(TN)\smile f^*\omega\bigr)[N]\;.$$
If~$n=m$, then~$\deg_{\hat A}$ is the usual degree.

Let~$D_N$ denote the untwisted Dirac operator on~$N$,
and let~$0\le\lambda_1(D_N^2)\le\lambda_2(D_N^2)\le\dots$
denote the eigenvalues of~$D_N^2$,
where each eigenvalue is repeated according to its multiplicity.

\begin{thm}\label{thm2}
  Let~$N$ be a closed Riemannian spin manifold and let~$k\in\N$.
  Then there is no Lipschitz map~$N\to S^{2m}$ of Lipschitz constant~$1$
  with
	$$\abs{\deg_{\hat A}f}>2^{m-1}\,(k-1)$$
  unless
	$$\lambda_k(D_N^2)\le\lambda_1(D_M)\;.$$
\end{thm}

If we replace the target manifold by a symmetric space~$G/H$
with~$\rk G=\rk H$,
then we can prove similar theorems where the degree
condition on~$f$ is replaced by
	$$\<[N],\hat A(TN)\smile f^*\alpha\>>C(k-1)\;,$$
with~$\alpha\in H^*(M;\Z)$ and~$C>0$ depending only on~$M$.
Alternatively,
there is no $1$-Lipschitz map~$N\to M$ with
	$$\abs{\deg_{\hat A}f}>C(k-1)$$
for some constant~$C$ depending on~$M$,
unless~$\lambda_k(D_N^2)\le c$ for some~$c$ depending only on~$M$,
with~$c>\lambda_1(D_M)$ in general.
This will be discussed in section~\ref{LipschitzSect}.
In the case~$\rk G-\rk H=1$ we need a $K$-theoretic condition on~$f$
instead of the cohomological $\hat A$-degree condition.
We may thus ask if even-dimensional spheres are the only
manifolds that admit optimal Vafa-Witten bounds for Lipschitz spin maps
of sufficiently large $\hat A$-degree.

\begin{rem}\label{rem-1}
The actual estimate and the index theoretic considerations
involved in the proof of Theorem~\ref{thm1} are very similar
to those used for the scalar curvature comparison result in~\cite{GS}.
Nevertheless,
we need a stronger metric condition
($\bar g\ge g$ on~$TM$, not just on~$\Lambda^2TM$).
This is due to the fact that the Vafa-Witten estimate
is related to Gromov's $K$-length,
whereas scalar curvature comparison is related to $K$-area,
see~\cite{Gromov}.
\end{rem}

\begin{rem}\label{rem0}
Note that Theorem~\ref{thm1} still holds if we replace~$\bar D$
by~$\bar D'=\bar D+A$,
where~$A$ is a symmetric endomorphism of the spinor bundle.
This is because in the proof,
we can then replace~$\bar D_i$ by~$\bar D_i'=\bar D_i+A\otimes\id$
for~$i=0$, $1$.
Then~$\bar D_1$ still has a kernel,
and of course~$\norm{\bar D'_0-\bar D'_1}=\norm{\bar D_1-\bar D_0}$.
\end{rem}

\begin{rem}\label{rem1}
Remark~\ref{rem0} in particular implies that on a symmetric space,
	$$\lambda_1(D+A)\le\lambda_1(D)$$
for all symmetric endomorphisms~$A$ of the spinor bundle.
In Section~\ref{BergerSect} we will exhibit an isotropy irreducible
quotient~$G/H=\SO(5)/\SO(3)$ of compact Lie groups with a normal metric,
for which the generalised estimate above does not hold.
This shows in particular that the methods of the present paper
do not readily generalise to normal homogeneous spaces.
\end{rem}

{\em Acknowlegdements.\/}
The rigidity statement in Theorem~\ref{thm1} was pointed out to us by
Listing.
An anonymous referee helped to make this article more readable.

\section{The Smallest Dirac Eigenvalue of a Symmetric Metric}\label{ParthSect}
Let~$M=G/H$ be symmetric quotient of compact Lie groups of equal rank
with Lie algebras~$\frh\subset\frg$.
We fix an Ad-invariant metric~$g$ on~$\frg$ and let~$\frp=\frh^\perp$.
The tangent bundle of~$M$ can be written as
	$$TM=G\times_H\frp\;,$$
where~$H$ acts on~$\frp$ by the restriction of the adjoint action~$\Ad_G$.
The scalar product~$g$ induces a symmetric Riemannian metric on~$M$
that will also be denoted by~$g$.
Note that if~$M$ is symmetric,
then the Levi-Civita connection on~$TM$ is precisely the reductive connection
on~$G\times_H\frp$.

Let~$\Sigma$ be a spinor module for the Clifford algebra of~$\frp$.
If we assume that~$M$ is $G$-spin,
then the $H$-representation of~$H$ on~$\frp$
induces an action~$\sigma\colon H\to\End\Sigma$.
The natural metric on~$\Sigma$ is $\sigma$-invariant,
so we obtain a $G$-equivariant metric~$g$ on~$SM$.
Equipped with this metric and the reductive connection,
the $G$-equivariant vector bundle
	$$SM=G\times_H\Sigma\to M$$
can be identified with the spinor bundle of~$M$,
and the $G$-equivariant Dirac operator~$D$
acts as an essentially selfadjoint operator on~$L^2(M;SM)$.

Let~$\hat G$ denote the set of equivalence classes of irreducible
$G$-representations.
We write~$\gamma\colon G\to\End V^\gamma$ for all~$\gamma\in\hat G$.
By Frobenius reciprocity and the Peter-Weyl theorem,
the $L^2$-sections of~$SM$ can be decomposed $G$-equivariantly
into a Hilbert sum
\begin{equation}\label{L2SMFormel}
  L^2(M;SM)
  =\overline{\bigoplus_{\gamma\in\hat G}
    V^\gamma\otimes\Hom_H(V^\gamma,\Sigma)}\;.
\end{equation}
The Dirac operator preserves this decomposition,
and we have
	$$D|_{V^\gamma\otimes\Hom_H(V^\gamma,\Sigma)}
	=\id_{V^\gamma}\otimes{}^{\gamma\!}D$$
with
\begin{equation}\label{DiracFormel}
  {}^{\gamma\!}D=\sum_{i=1}^n\gamma^*_{e_i}\otimes c_i
  \quad\in\quad\End(\Hom_H(V^\gamma,\Sigma))\;.
\end{equation}
Here~$e_1$, \dots, $e_n$ is a $g$-orthonormal base of~$\frp$
and~$c_i$ denotes Clifford multiplication by~$e_i$.

Let~$c_G^\gamma$ denote the Casimir operator of the $G$-action~$\gamma$
with respect to the metric~$g$,
which acts as a scalar on~$V^\gamma$.
If~$M$ is symmetric,
then the Casimir operator~$c_H^\sigma$ of~$\sigma$ also acts as a scalar,
even though~$\sigma$ is in general not irreducible.
By Parthasarathy's formula~\cite{Parth},
	$${}^{\gamma\!}D^2=c_G^{\gamma^*}+c_H^\sigma\;.$$

\begin{Proposition}\label{SymmEVProp}
The smallest eigenvalue of~$D^2$
is given by
	$$\lambda_1(D^2)
	=\min\bigl\{\,c_G^{\gamma^*}+c_H^\sigma
	  \bigm|\gamma\in\hat G\text{ with }
	    \Hom_H(V^{\gamma},\Sigma)\ne 0\,\bigr\}\;.\quad\qed$$
\end{Proposition}

A more explicit formula for the first eigenvalue in the case~$\rk G=\rk H$
has recently been given in~\cite{Milhorat}.

\section{The Vafa-Witten Estimate}\label{EstSect}
Let~$M=G/H$ and~$g$ as before
and consider an arbitrary Riemannian metric~$\bar g$ on~$M$ with~$\bar g>g$.
The corresponding Dirac operator will be denoted~$\bar D$.
Let~$\lambda_1(\bar D^2)$ denote the smallest eigenvalue of~$\bar D^2$.
As mentioned in the introduction,
we will estimate~$\lambda_1(\bar D^2)$
by comparing the related operator~$\bar D_0=\bar D\otimes\id_{\C^N}$
to an operator~$\bar D_1$ acting on the same space of sections.

Assume that we are given vector bundle~$W\subset V=M\times\C^N\to M$
such that the twisted Dirac operator on~$SM\otimes W$ has nonvanishing
index;
these bundles will be constructed in steps below
and in sections~\ref{eqr} and~\ref{oddsect}.
Let~$\nabla^0$ be the trivial connection on~$V$,
and let~$\nabla^1$ be another connection for which~$W\subset V$ is a parallel
subbundle.
Let~$\bar D_0$, $\bar D_1$ denote the corresponding twisted Dirac operators
on~$SM\otimes\C^N$ for the metric~$\bar g$.
Then~$\bar D_0$ is just the direct sum of~$N$ copies of~$\bar D$,
whereas~$\bar D_1$ has nontrivial kernel.
By a Rayleigh quotient argument,
thus
\begin{equation}\label{VaWiEst}
  \lambda_1(\bar D^2)
  =\lambda_1(\bar D_0^2)
  \le\frac{\norm{(\bar D_1-\bar D_0)s}_{L^2}^2}{\norm s_{L^2}^2}
  \le\sup_{p\in M}\norm{(\bar D_1-\bar D_0)_p^2}_{\mathrm{op}}
\end{equation}
for some~$0\ne s\in\ker(\bar D_1)$.
Here~$\norm{\,\cdot\,}_{L^2}$ denotes the $L^2$-norm of sections,
where\-as~$\norm{\,\cdot\,}_{\mathrm{op}}$ denotes the pointwise operator norm
of an endomorphism.
The operator~$\bar D_1-\bar D_0$ is of order zero,
so~$\norm{(\bar D_1-\bar D_0)_p^2}_{\mathrm{op}}$
is well-defined and can be estimated pointwise on~$M$.

In our current situation,
let~$\gamma\in\hat G$ be a $G$-representation such that
\begin{equation}\label{GammaChoice}
  \Hom_H(V^\gamma,\sigma)\ne 0
  \qquad\text{and}\qquad
  \lambda_1(D^2)=c_G^{\gamma^*}+c_H^\sigma\;,
\end{equation}
cf.\ Proposition~\ref{SymmEVProp}.
Let~$\gamma^*$ denote the dual representation on~$V^{\gamma^*}=(V^\gamma)^*$
and let
	$$V=G\times_H(V^{\gamma^*}|_H)$$
denote the corresponding vector bundle over~$G$,
then~$V$ is trivialized by the map
	$$V\to M\times V^{\gamma^*}
		\qquad\text{with}\qquad
	[g,v]\mapsto (gH,\gamma^*_gv)\;.$$
To a map~$v\colon M\to V^\gamma$ corresponds the section
	$$gH\mapsto[g,(\gamma^*)^{-1}_g v(g)]\in V_{gH}\;.$$
Now let~$\nabla^0$ denote the trivial connection
and let~$\nabla^1$ denote the reductive connection on~$V$.
Then~$\overline D_0$, $\overline D_1$ are the corresponding
Dirac operators on~$SM\otimes V$.

If we apply the reductive connection~$\nabla^1$ on~$V$ to a section that
is constant in the given trivialisation,
then
	$$\nabla^1_{[g,X]}v
	=\bigl[g,\gamma^*_X(\gamma^*)^{-1}_gv\bigr]=(gH,\gamma_{*\Ad_gX}v)$$
where~$X\in\frp=\frh^\perp$ and
	$$[g,X]=\frac\del{\del t}\,\bigl(ge^{-tX}H\bigr)\in T_{gH}M\;.$$
Thus we can identify~$T_pM$ with~$\frp$
and choose a $g$-orthonormal frame~$f_1$, \dots, $f_m$ of~$\frh$.
We recall that the action~$\sigma_*$ of~$\frh$ on~$SM$ can
be described in terms of Clifford multiplication and Lie brackets by
	$$\sigma_{*f_k}=\frac14\sum_{i,j=1}^n\<[e_i,e_j],f_k\>\,c_ic_j\;.$$

Let~$\bar e_1,\dots,\bar e_n$ be a $\bar g$-orthonormal base of~$T_pM$.
We may assume that there exist~$0<\mu_1$, \dots, $\mu_n\le 1$
such that~$\bar e_i=\mu_ie_i$ for all~$i\in\{1,\dots,n\}$,
where~$e_1$, \dots, $e_n$ is an orthonormal base with respect to~$g$.
We have to compute the operator norm of the operator
	$$C=\bar D_1-\bar D_0
	=\sum_{i=1}^nc_i\gamma^*_{\bar e_i}
	=\sum_{i=1}^n\mu_ic_i\gamma^*_{e_i}\;,$$
which is the difference of two Dirac operators on~$SM\otimes V$
with respect to the metric~$\bar g$.
We note that this formula looks similar to~\eqref{DiracFormel} above.
This is a special feature of symmetric spaces,
which does not even generalise to normal homogeneous spaces,
see Section~\ref{BergerSect} below.

Because~$C$ is selfadjoint,
we have~$\norm{Cv}^2=\<C^2v,v\>$,
and it suffices to estimate the eigenvalues of~$C^2$.
We follow the proof of Parthasarathy's formula~\cite{Parth}.
Using~$[\frp,\frp]\subset\frh$
we compute
\begin{align*}
  C^2
  &=-\sum_{i=1}^n\mu_i^2(\gamma^*_{e_i})^2
    +\frac12\,\sum_{i,j=1}^n\mu_i\mu_j\,c_ic_j\,\gamma^*_{[e_i,e_j]}\\
  &=-\sum_{i=1}^n\mu_i^2(\gamma_{e_i}^*)^2
    +\sum_{k=1}^m\biggl(\gamma^*_{f_k}
      +\frac14\sum_{i,j=1}^n\mu_i\mu_j
	\,\<[e_i,e_j],f_k\>\,c_ic_j\biggr)^2\\
  &\qquad
      -\sum_{k=1}^m(\gamma^*_{f_k})^2
      -\frac1{16}\sum_{k=1}^m\biggl(\sum_{i,j}\mu_i\mu_j
         \,\<[e_i,e_j],f_k\>\,c_ic_j\biggr)^2\\
  &\le c_G^{\gamma^*}
      -\frac1{16}\sum_{i,j,k,l=1}^n\mu_i\mu_j\mu_k\mu_l
         \,\<[e_i,e_j],[e_k,e_l]\>\,c_ic_jc_kc_l\;.
\end{align*}

A term similar to the last one on the right hand side has already been
estimated in~\cite{GS}, equation~(1.11).
Because~$[\frp,\frp]\subset\frh$, we have
\begin{multline}\label{ParthEst}
  -\frac1{16}\sum_{i,j,k,l=1}^n\mu_i\mu_j\mu_k\mu_l
         \,\<[e_i,e_j],[e_k,e_l]\>\,c_ic_jc_kc_l
  =\frac18\sum_{i,j=1}^n\mu_i^2\mu_j^2\,\norm{[e_i,e_j]}^2\\
  \le\frac18\sum_{i,j=1}^n\norm{[e_i,e_j]}^2
  =-\sum_{k=1}^m\biggl(\sum_{i,j=1}^n\<[e_i,e_j],f_k\>\,c_ic_j\biggr)^2
  =c_H^\sigma\;.
\end{multline}

Combining the calculations above,
we have the estimate
\begin{equation}\label{GammaEst}
  C^2=(\bar D_1-\bar D_0)^2
  \le c_G^{\gamma^*}+c_H^\sigma
  =\lambda_1(D^2)\;.
\end{equation}

\section{Rigidity}\label{RigSect}

We will now give a short of proof due to M. Listing
that the equality~$\lambda_1(\bar D^2)=\lambda_1(D^2)$
in Theorem~\ref{thm1} implies that~$\bar g=g$.

Let~$0\ne s\in\ker(\bar D_1)$ as in~\eqref{VaWiEst}.
If~$\lambda_1(\bar D^2)=\lambda_1(D^2)$,
we conclude that
	$$\norm{(\bar D_1-\bar D_0)s}_{L^2}^2
	=\lambda_1(D^2)\,\norm s^2_{L^2}\;.$$
by~\eqref{VaWiEst} and~\eqref{GammaEst}.
This implies in particular
that~$\norm{(\bar D_1-\bar D_0)_p^2}_{\mathrm{op}}=\lambda_1(D^2)$
holds for all~$p\in\operatorname{supp}(s)$.
By~\cite{Baer}, we know that~$s$ is nonzero on a dense open subset of~$M$,
so we have
\begin{equation}\label{ParthEq}
  \norm{(\bar D_1-\bar D_0)^2_p}_{\mathrm{op}}=\lambda_1(D^2)
\end{equation}
for all~$p\in M$,
and we must have equality in~\eqref{ParthEst}.

Since~$M$ is of compact type, $\mathfrak g$ has trivial center.
In particular,
for each~$1\le i\le n$ there exists~$1\le j\le n$ such that~$[e_i,e_j]\ne 0$.
Thus~\eqref{ParthEq} implies by~\eqref{ParthEst} that~$\mu_1=\dots=\mu_n=1$.
The rigidity statement follows once we have shown that~$\ker(\bar D_1)\ne 0$.

\section{The Equal Rank Case}\label{eqr}
We prove Theorem~\ref{thm1} for~$\rk G=\rk H$.
It remains to show that the operator~$\bar D_1$ has a kernel.

\begin{Proposition}[\cite{Parth}, \cite{Goette}]\label{ParthProp}
Let~$M=G/K$ be a symmetric space with~$\rk G=\rk K+k$.
Then the complex spinor bundle~$\Sigma$
is locally induced by a representation~$\sigma$
of the Lie algebra~$\mathfrak k$ of~$K$,
which splits as
        $$\sigma=2^{\textstyle\left[\frac{k}{2}\right]}_{\phantom{x}}
                        \,\bigoplus_{i=1}^q\sigma_i\;,$$
where~$\sigma_1$, \dots, $\sigma_q$ are certain pairwise non-isomorphic
irreducible complex representations of~$\mathfrak k$.
\end{Proposition}

By our choice of~$\gamma$ in~\eqref{GammaChoice},
we have~$\Hom_H(V^\gamma,\Sigma)\ne 0$.
By Schur's lemma,
the $H$-representations~$\gamma|_H$ and~$\sigma$
contain a common irreducible subrepresentation,
say~$\sigma_1$ acting on~$\Sigma_1$.
Consider the vector bundle
	$$W=G\times_H\Sigma_1^*\;,$$
then~$W$ is a parallel subbundle of~$V$
with respect to the reductive connection~$\nabla^1$.
Let~$D_1$ be the Dirac operator acting on~$SM\otimes W$.

As above, we now have
	$$L^2(M;SM\otimes W)
	=\bigoplus_{\gamma\in\hat G}
		V^\gamma\otimes\Hom_H(V^\gamma,\Sigma\otimes\Sigma_1^*)\;.$$
Again by Parthasarathy's formula (\cite{Parth}, \cite{Goette}),
the operator~$D_1^2$ acts on~$V^\gamma,\Sigma\otimes\Sigma_1$
as
	$${}^{\gamma\!}D_1^2=c_G^\gamma+c_H^\sigma-c_H^\sigma=c_G^\gamma\;.$$
Because~$c_G^\gamma=0$ iff~$\gamma$ is the trivial representation,
the kernel of~$D_1^2$ is precisely
	$$\ker D_1=\ker D_1^2=
	\Hom_H(\C,\Sigma\otimes\Sigma_1^*)\cong\Hom_H(\Sigma_1,\Sigma)
	=\Hom_H(\Sigma_1,\Sigma_1)$$
by Schur's Lemma.
In particular~$\dim\ker D_1=1$,
whence~$\ind D_1=\pm 1\ne 0$.

To complete the proof of Theorem~\ref{thm1} in this case,
we note that the Dirac operator on~$\overline{SM}\otimes W$
has nonzero index,
and hence~$\bar D_1$ has a kernel.
The claim now follows from~\eqref{VaWiEst}
and~\eqref{GammaEst}.

\section{Mod-2-Indices and the Case \texorpdfstring{$\rk H=\rk G-1$}{of rank difference one}}\label{oddsect}
We now prove Theorem~\ref{thm1} for~$\rk G=\rk H+1$.
We proceed similar as in~\cite{GS}, section~2.c.

Assume first that~$n=\dim M\equiv 1$ mod~$8$.
In this case,
there exists a real vector bundle~$S_\R M$ such that~$SM=S_\R M\otimes_\R \C$,
and the even part of the real Clifford algebra still acts on~$S_\R M$.
The bundle~$S_\R M$ is induced by real representation~$\sigma_\R$ of~$H$.
Let~$\sigma_{\R,i}$ denote its irreducible components,
such that~$\sigma_i=\sigma_{\R,i}\otimes_\R \C$,
and let~$W_\R=G\times_H\Sigma_{\R,1}^*$.

Let
	$$\omega_\R=c_1\cdots c_n$$
denote the real Clifford volume element,
then~$\omega_\R$ is parallel, anti-selfadjoint,
and commutes with~$D\omega_\R^*=\omega_\R$.
The operator
	$$D_{\R,1}=\omega_\R\,D_1$$
is anti-selfadjoint and
has coefficients in the even part of the real Clifford algebra,
so it acts on~$S_\R M$.
By the same reasoning as above,
	$$\ker(\omega_\R\,D_{\R,1})
	=\Hom_{\R,H}(\C,\Sigma_{\R,1}\otimes_\R \Sigma_{\R,1}^*)
	=\Hom_{\R,H}(\Sigma_{\R,1},\Sigma_{\R,1})$$
is one-dimensional.

With respect to the metric~$\bar g$,
we similarly construct the operator~$\bar D_{\R,1}=\bar\omega_\R\,\bar D_1$
acting on~$\bar S_\R M\otimes_\R W_\R$.
Its restriction to~$\overline S_\R M\otimes_\R V_\R$ can be deformed
into~$D_{\R,1}$ through a family of elliptic, formally
anti-selfadjoint operators.
Since the parity of the dimension of the kernel is preserved by such a
deformation (see~\cite{LM}),
the operator~$\bar\omega_\R\,\bar D_{\R,1}$ has nontrivial kernel.
We regard~$V\to M$ as a real vector bundle by forgetting the complex structure,
then~$W_\R$ is again a parallel subbundle of~$V$ with respect to~$\nabla^1$,
so~$\bar D_1$ has a nontrivial kernel, too.
Once again,
our theorem follows from~\eqref{VaWiEst} and~\eqref{GammaEst}.

Finally,
if~$\rk G-\rk H=1$ and~$\dim M\not\equiv 1$ mod~$8$,
we consider~$N=M\times S^{2m}$ with~$\dim N\equiv 1$ mod~$8$.
We fix a round metric~$g_0$ on~$S^{2m}$ and regard~$N$
with the metrics~$g\oplus g_0$ and~$\bar g\oplus g_0$.
Then our main theorem holds for~$N$,
and the result for~$M$ follows because
	$$\lambda_1(\bar D_M^2)
	=\lambda_1(\bar D_N^2)-\lambda_1(D_{S^{2m}}^2)
	\le\lambda_1(D_N^2)-\lambda_1(D_{S^{2m}}^2)
	=\lambda_1(D_N^2)\;.\quad\qed$$

\section{Examples}\label{ExpSect}
We consider spheres and projective spaces.

\begin{Example}\label{SnEx}
For even-dimensional spheres~$S^{2m}=\Spin(2m+1)/\Spin(2m)$,
the restriction of the $2^m$-dimensional spinor representation~$\sigma$
of~$\Spin(2m+1)$ to~$H$ splits into two irreducible representations
	$$\sigma|_{\Spin(2m)}=\sigma^+\oplus\sigma^-\;.$$
In~\cite{herz}, the bundles
	$$V=SM=G\times_H\Sigma\quad\supset\quad W=S^+M=G\times_H\Sigma^+$$
were used to prove the optimality of the Vafa-Witten estimate.
Our method is a direct generalisation of this approach
to other symmetric spaces.
\end{Example}

\begin{Example}\label{CPnEx}
The complex projective space~$\C P^n$ is spin iff~$n$ is odd,
so we regard the odd-dimensional projective space~$M=\C P^{2m-1}=G/H$
with~$G=\SU(2m)$ and~$H=\U(2m-1)$.
The Dirac spectrum has been computed in~\cite{CFG} and~\cite{Seeger}.
To summarize,
the spinor bundle splits as
	$$SM=\bigoplus_{q=0}^{2m-1}\Lambda^{0,q}T^*M\otimes\tau^m\;,$$
where~$\Lambda^{0,q}T^*M$ denotes the bundle of anti-holomorphic
differential forms of degree~$q$,
and~$\tau$ denotes the tautological bundle.
Note that~$\tau^m$ is a square root of the canonical bundle.

The smallest eigenvalue of~$D^2$ on $\Lambda^{0,q}T^*M\otimes\tau^m$
with respect to the Fubini-Study metric is given by
\begin{equation}\label{CPnEV}
  \lambda_1\bigl(D^2|_{\Lambda^{0,q}T^*M\otimes\tau^m}\bigr)
  =\begin{cases}
    8m^2-4m(q+1)&\text{for~$q<m$, and}\\
    8m^2-4m(2m-q)&\text{for }q\ge m\;.
  \end{cases}
\end{equation}
Thus the lowest eigenvalue is attained in the middle degrees~$q=m-1$, $m$,
and is given by~$4m^2=(n+1)^2$.
In contrast, Davaux and Min-Oo used
the bundle~$\tau^m$ for~$q=0$ in~\cite{DM},
which explains their larger upper bound~$8m^2-4m=2n(n+1)$.

Let us give a geometric description
of the bundles~$\Lambda^{0,q}T^*M\otimes\tau^m$ for~$q=m-1$, $m$.
We have
\begin{align*}
  T\C P^{2m-1}
  &=\bigl\{\,[x,v]\in S^{2m-1}\otimes\C^{2m}\bigm|x\ne0,x\perp v\,\}\\
  \text{with}\qquad
       [z,v]&=\bigl\{\,(zx,zv)\bigm|z\in S^1\,\bigr\}\;.
\end{align*}
For~$\alpha\in\Lambda^{0,m}T^*M$,
we note that
\begin{align*}\alpha(zx,zv_1,\dots,zv_{m-1})
  &=\bar z^m\,\alpha(x,v_1,\dots,v_{m-1})\\
  \text{and}\qquad
  \alpha(zv_1,\dots,zv_m)
  &=\bar z^m\,\alpha(v_1,\dots,v_m)\;.
\end{align*}
Thus inserting representatives of tangent vectors into $m$-forms gives
a natural trivialisation
	$$\Phi=\Phi^{m-1}\oplus\Phi^m
	\colon\C P^{2m-1}\times\Lambda^{0,m}\C^{2m}
		\longrightarrow
	(\Lambda^{0,m-1}T^*M\oplus\Lambda^{0,m}T^*M)\otimes\tau^m\;.$$
\end{Example}

\section{Lipschitz maps of High \texorpdfstring{$\hat A$}{A}-Degree}\label{LipschitzSect}
We give a proof of Theorem~\ref{thm2}.
Let~$N$ be a closed Riemannian spin manifolds and~$f\colon N\to M$
a Lipschitz map of Lipschitz constant~$1$.
Then~$f$ can be approximated by smooth $(1+\eps)$-Lipschitz maps
for all~$\eps>0$.
We may therefore assume that~$f$ is smooth.

As in~\cite{baum}, \cite{herz}, consider
the spinor bundles~$\Sigma^\pm\to S^{2n}$
equipped with the Dirac connection~$\nabla^1$
with respect to the round metric~$g$.
Let~$\nabla^0$ denote the trivial connection on~$\Sigma^+\oplus\Sigma^-
\cong S^{2n}\otimes\Sigma$,
where~$\Sigma$ is the spinor module of~$\R^{2n+1}$.

The Chern character of the spinor bundles is given by
	$$\ch(\Sigma^\pm)=2^{m-1}\pm\omega\;,$$
where~$\omega\in H^{2m}(S^{2m};\Z)$ is a generator.
Let~$D_{N,1,\pm}$ be the Dirac operator on~$N$ twisted by~$f^*\Sigma^\pm$,
and let~$D_{N,1}$ denote their direct sum.
By the Atiyah-Singer index theorem,
we have
	$$\ind(D_{N,1,\pm})=\bigl(\hat A(TN)\smile f^*\ch(\Sigma^\pm)\bigr)[N]
	=2^{m-1}\hat A(TN)[N]\pm\deg_{\hat A}f\;.$$
By combining both bundles,
it is now easy to see that
\begin{equation}\label{DimKer}
  \dim\ker(D_{N,1})
  \ge\max\bigl(\bigl|2^m\hat A(TN)[N]\bigr|,2\abs{\deg_{\hat A}f}\bigr)
  >2^m(k-1)
\end{equation}
if~$\abs{\deg_{\hat A}f}>2^{m-1}(k-1)$.

Let~$p\in N$,
then there exists an orthonormal frame~$\bar e_1,\dots,\bar e_n$
of~$T_pN$ and an orthonormal frame~$e_1, \dots, e_{2m}$ of~$T_{f(p)}S^{2m}$
such that
	$$d_pf(\bar e_k)
	=\begin{cases}
		\mu_k\,e_k&\text{if~$k\le2m$, and}\\
		0&\text{otherwise.}
	\end{cases}$$
If~$f$ is~$(1+\eps)$ Lipschitz,
then~$\mu_1$, \dots, $\mu_{2m}\le1+\eps$.

Let~$D_{N,0}$ denote the Dirac operator on~$N$ twisted by~$f^*\Sigma$,
but with respect to the trivial connection~$\nabla^0$.
Then
	$$D_{N,1}-D_{N,0}
	=\sum_{i=1}^{2m}\bar c_i\,\gamma_{*\mu_ie_i}\;.$$
Now a similar computation as in Section~\ref{EstSect} gives
\begin{equation}\label{MapVaWi}
  \norm{D_{N,1}-D_{N,0}}^2\le (1+\eps)\,(c_G^\gamma+c_H^\sigma)
  =(1+\eps)\,\lambda_1(D^2)\;.
\end{equation}

Combining~\eqref{DimKer} and~\eqref{MapVaWi} as before,
we conclude at least~$2^m(k-1)$ eigenvalues of~$D_{N,0}$
are not larger than~$(1+\eps)\,\lambda_1(D^2)$.
Because~$D_{N,0}$ consists of~$2^m$ copies of the Dirac operator~$D_N$,
we conclude that at least~$k$ eigenvalues of~$D_N^2$
(counted with multiplicities) are not larger than~$(1+\eps)\,\lambda_1(D^2)$.
Since we can choose~$\eps>0$ arbitrarily small,
Theorem~2 is proved.\qquad\qed

\medskip

We remarked in the Introduction that Theorem~\ref{thm2} does not hold
unchanged for other symmetric spaces.
Regard~$M=\C P^{2m-1}$ as in Example~\ref{CPnEx}.
By~\eqref{CPnEV},
the Vafa-Witten will be optimal only if we choose
one of the bundles~$W=\Lambda^{0,m-1}T^*M\otimes\tau^m$
or~$W^*=\Lambda^{0,m}T^*M\otimes\tau^m$ as twist bundle.

To determine the Chern characters of these bundles,
recall that
	$$T'\C P^{2m-1}\oplus\C\cong2m\,\tau^{-1}\;,$$
and thus in $K$-theory,
	$$[\Lambda^{0,q}T^*M]
	=\sum_{i=0}^q(-1)^{q-i}\,[\Lambda^i(2m\,\tau^{-1})]
	=\sum_{i=0}^q(-1)^{q-i}\,\begin{pmatrix}2m\\i\end{pmatrix}
		\,[\tau^{-i}]\;.$$
We know that~$a=c_1(\tau)$ is a generator of~$H^2(\C P^{2m-1})$,
and by the above, we find
	$$\ch(W)
	=\sum_{i=0}^{m-1}(-1)^{m-1-i}\,\begin{pmatrix}2m\\i\end{pmatrix}
		\,e^{(m-i)a}\;.$$
Already for~$\C P^3$, the explicit classes are
	$$\ch(W)=3+2a-\frac23\,a^3
		\qquad\text{and}\qquad
	\ch(W^*)=3-2a+\frac23\,a^3\;.$$

On the other hand,
we have seen that
	$$W\oplus W^*\cong
	\C^{\textstyle\bigl(\begin{smallmatrix}2m\\
		m\end{smallmatrix}\bigr)}\;.$$
Proceeding as above,
we can prove the following result.

\begin{Proposition}\label{CPnProp1}
  Let~$N$ be a closed Riemannian spin manifold and let~$k\in\N$.
  Then there is no Lipschitz map~$N\to\C P^{2m-1}$ of Lipschitz constant~$1$
  with
	$$\biggl|\hat A(TN)\smile f^*\sum_{i=0}^{m-1}(-1)^i
		\,\begin{pmatrix}2m\\i\end{pmatrix}
		\,\bigl(e^{(m-i)a}-1\bigr)\biggr|
	>\begin{pmatrix}2m-1\\m-1\end{pmatrix}\,(k-1)\;,$$
  where~$a$ is a generator of~$H^2(\C P^{2m-1})$,
  unless
	$$\lambda_k(D_N^2)\le\lambda_1(D_M^2)\;.\quad\qed$$
\end{Proposition}

\section{A Homogeneous Counterexample}\label{BergerSect}
We now explain Remark~\ref{rem1}.
In general,
the spectrum of a Dirac operator on a homogeneous space is hard to compute.
On the Berger space~$M=\SO(5)/\SO(3)$ however,
most of the relevant geometric structure can be described using
octonions, see~\cite{GKS}.
We will see that a straightforward adaptation of the arguments
of Sections~\ref{ParthSect}--\ref{oddsect} to the Berger space is not possible.

We embed~$H=\SO(3)$ in~$G=\SO(5)$ via the irreducible $\SO(3)$-representation
of dimension~$5$.
On~$\frg=\so(5)$,
we take the scalar product~$\<A,B\>=-\frac12\,\tr(AB)$.
As above, let~$\frp$ denote the orthogonal complement of~$\frh=\so(3)$.
Then~$\dim\frp=7$, and the isotropy action of~$\SO(3)$ on~$\frp$ factors as
\begin{equation}\label{BergerG2}
  \SO(3)\longrightarrow\G_2\longrightarrow\Spin(7)
  \longrightarrow\SO(7)\;.
\end{equation}
Let~$[\cdot,\cdot]_{\frp}$ denote the projection of the Lie bracket
in~$\frg$.
Let~$\I\subset\O$ denote the imaginary octonions,
let~$*$ denote octonion multiplication,
and let~$*_\I$ denote octonion multiplication followed by the projection
onto~$\I$.

\begin{Lemma}[\cite{GKS}]\label{GKS1}
With a suitable isometric, $\G_2$-equivariant identification of~$\frp$
with the imaginary octonions,
one has
	$$[v,w]_\frp=\frac1{\sqrt 5}\,v\mathbin{*_\I}w
	\qquad\text{for all}\qquad v, w\in\frp\;.\quad\qed$$
\end{Lemma}

Let~$\Sigma_\R$ again denote the real spinor module of~$\frp$,
then~$\G_2$ also acts on~$\Sigma_\R$ by~\eqref{BergerG2}.

\begin{Lemma}[\cite{GKS}]\label{GKS2}
We identify $\frp\cong\I$ as in Lemma~\ref{GKS1}. With respect to a suitable
isometric, $\G_2$-equivariant identification $\Sigma_\R\cong\O$ and a suitable
orientation of~$\frp$,
Clifford multiplication~$\frp\times\Sigma_\R\to\Sigma_\R$
equals Cayley multiplication $\I\times\O\to\O$ from the right.\qed
\end{Lemma}

We fix an orthonormal base~$e_1$, \dots, $e_7$ of~$\frp\cong\I$ such that
	$$e_i*e_{i+1}=e_{i+3}\;,$$
where indices are taken mod~$7$.
Clifford multiplication with~$e_i$ will be abbreviated as~$c_i$.
As in~\cite{Goette} and~\cite{GKS},
let us introduce the notation
\begin{align*}
  c_{ijk}&=\<[e_i,e_j]_\frp,e_k\>
  =-\frac1{\sqrt 5}\Re\bigl((e_i*e_j)*e_k\bigr)\\
  \text{and}\qquad
  \adsl_{\frp,X}&=\frac14\,\sum_{i,j=1}^7\<[X,e_i],e_j\>\,c_ic_j\;.
\end{align*}

Let us define
	$$A=\frac13\sum_{i=1}^7c_i\adsl_{\frp,e_i}
	=\frac1{12}\sum_{i,j,k=1}^7c_{ijk}\,c_ic_jc_k\;.$$
It is easy to see that~$A$ is a symmetric, $\G_2$-invariant
endomorphism of~$\Sigma_\R$.
Thus it induces a $G$-invariant symmetric endomorphism of~$S_\R M$.
One may compute that
	$$A=\frac1{2\sqrt 5}\begin{pmatrix}7&0\\0&-1\end{pmatrix}
	\colon\R\oplus\I\to\R\oplus\I\;.$$
Let~$D$ denote the real Dirac operator acting on~$\Gamma(S_\R M)$.
We want to consider the family of operators
\begin{equation}\label{DlambdaDef}
  D^\lambda
  =D+\biggl(3\lambda-\frac32\biggr)\,A\;.
\end{equation}
The operator~$D^{\frac12}$ is the Riemannian Dirac operator on~$M$,
whereas~$\tilde D=D^{\frac13}$ is called the {\em reductive\/}
or {\em cubic\/} Dirac operator.

\begin{Proposition}\label{BergerProp}
  For~$\lambda$ sufficiently close to~$\frac12$,
  the smallest eigenvalue of~$(D^\lambda)^2$
  is~$\frac{441}{20}\,\lambda^2$.
\end{Proposition}

In particular there exists~$\eps>0$ sufficiently small such that
	$$\lambda_1\bigl((D^{\frac12+\eps})^2\bigr)>\lambda_1(D^2)\;.$$
As stated in Remark~\ref{rem1},
this implies that the Vafa-Witten technique cannot be used
to prove sharp upper bounds for the first Dirac eigenvalue.

For certain questions,
the reductive Dirac operator~$\tilde D$ on a normal homogeneous space
plays the same role as the geometric Dirac operator on a symmetric space,
but here clearly~$\lambda_1(\tilde D^2)<\lambda_1(D^2)$,
so with respect to optimal Vafa-Witten estimates,
the reductive Dirac operator is even worse than the geometric Dirac operator.

\begin{proof}
  The space of sections~$L^2(M;S_\R M)$ can be decomposed
  as in~\eqref{L2SMFormel}.
  As in~\cite{Goette},
  we see that~$D^\lambda$ acts on the isotypical
  components~$V^\gamma\otimes\Hom_H(V^\gamma,\Sigma_\R)$
  by~$\id_{V^\gamma}\otimes{}^{\gamma\!}D^\lambda$,
  where
	$${}^{\gamma\!}D^\lambda
	=\sum_{i=1}^7\Bigl(\gamma_{e_i}^*\otimes c_i
  		+\lambda\,\id_{V^\gamma}\otimes c_i\adsl_{\frp,e_i}\Bigr)
	=\sum_{i=1}^7\gamma_{e_i}^*\otimes c_i
		+\frac\lambda4\sum_{i,j,k=1}^7c_{ijk}\,c_ic_jc_k\;.$$

  Let us write~$\mu=3\lambda-1$,
  then
	$$D^\lambda=\tilde D+\mu\,\id_{V^\gamma}\otimes A\;.$$
  The square of~$\tilde D$ has been computed in~\cite{Goette} as
	$$({}^{\gamma\!}\tilde D)^2
	=\norm{\gamma+\rho_G}^2-\norm{\rho_H}^2\;.$$
  If we use integers~$p\ge q\ge 0$ to index the irreducible representations
  of~$\SO(5)$,
  then by~\cite{GKS},
  we have the explicit formula
	$$({}^{\gamma_{p,q}\!}\tilde D)^2
	=\norm{\gamma_{p,q}+\rho_G}^2-\norm{\rho_H}^2
	=p^2+3p+q^2+q+\frac{49}{20}
		\quad\text{on }\Hom_H(V^{\gamma_{p,q}},\Sigma_\R)\;.$$

  Let us now analyze all irreducible representations of~$\SO(5)$.
  For the trivial representation, clearly
	$${}^{\gamma_{0,0}\!}D^\lambda
	=3\lambda\,{}^{\gamma_{0,0}\!}\tilde D
	=3\lambda\,A|_{\R}
	=\frac{21}{2\sqrt 5}\,\lambda\;,
		\qquad\text{so}\qquad
	{}^{\gamma_{0,0}\!}D={}^{\gamma_{0,0}\!}D^{\frac12}
	=\frac{21}{4\sqrt 5}$$
  The standard representation~$\gamma_{1,0}$ contains no $\SO(3)$-irreducible
  subrepresentation isomorphic to~$\R$ or~$\I$,
  so it does not contribute to the spectrum of~$D^\lambda$.
  The representation~$\gamma_{1,1}=\Lambda^2\gamma_{1,0}$ contains
  no trivial $\SO(3)$-subrepresentation,
  but an $\SO(3)$-subrepresentation isomorphic to~$\I$,
  and we find that
	$$\abs{\lambda_1({}^{\gamma_{1,1}\!}D)}
	\ge\bigl|{}^{\gamma_{1,1}\!}\tilde D\bigr|
		-\frac12\,\bigl|A|_\I\bigr|
	=\frac{13}{2\sqrt 5}-\frac1{4\sqrt 5}
	>\frac{21}{4\sqrt 5}=\abs{{}^{\gamma_{0,0}\!}D}\;.$$
  For all other representations, we have~$p>2$, thus
	$$\abs{\lambda_1({}^{\gamma_{p,q}\!}D)}
	\ge\bigl|{}^{\gamma_{2,0}\!}\tilde D\bigr|
		-\frac12\,\abs{A}
	=\frac{\sqrt{249}}{2\sqrt 5}-\frac7{4\sqrt 5}
	>\frac{21}{4\sqrt 5}=\abs{{}^{\gamma_{0,0}\!}D}\;.$$
  In particular,
  the smallest eigenvalue of~$D^\lambda$ equals~$\frac{21}{2\sqrt 5}\,\lambda$
  for~$\lambda$ sufficiently close to~$\frac12$,
  and is attained by a $G$-invariant spinor.
\end{proof}

\end{document}